\documentclass[12pt]{amsart}
\usepackage{amsmath, amsthm, amscd, amsfonts}

\newtheorem{theorem}{Theorem}[section]

\newtheorem{proposition}[theorem]{Proposition}
\newtheorem{corollary}[theorem]{Corollary}
\theoremstyle{definition}

\theoremstyle{remark}
\newtheorem{remark}[theorem]{Remark}
\numberwithin{equation}{section}
\begin{document}
\title{Satbility of Ternary Homomorphisms via Generalized Jensen Equation}
\author{M. S. Moslehian}
\address{Mohammad Sal Moslehian\newline Department of Mathematics, Ferdowsi University, P. O. Box 1159, Mashhad 91775, Iran}
\email{moslehian@ferdowsi.um.ac.ir}
\author{L. Sz\' ekelyhidi}
\address{L\' aszl\' o Sz\' ekelyhidi \newline Institute of Mathematics, Lajos Kossuth University, 4010 Debrecen,
Hungary} \email{szekely@math.klte.hu}
\begin{abstract}
In this paper, we establish the generalized Hyers--Ulam--Rassias
stability of homomorphisms between ternary algebras associted to
the generalized Jensen functional equation $r f(\frac{sx+ty}{r}) =
s f(x) + t f(y)$.
\end{abstract}
\maketitle

\section {Introduction}
A ternary (associative) algebra $({\mathcal A}, [ \; ])$ is a
linear space ${\mathcal A}$ over a scalar field ${\mathbb F} =
{\mathbb R} ~\text{or}~ {\mathbb C}$ equipped with a linear
mapping, the so-called ternary product, $[ \; ] : {\mathcal A}
\times {\mathcal A} \times {\mathcal A} \to {\mathcal A}$ such
that $[[abc]de] = [a[bcd]e] = [ab[cde]]$ for all $a, b, c, d, e\in
{\mathcal A}$. This notion is a natural generalization of the
binary case. Indeed if $({\mathcal A}, \odot)$ is a usual
(binary) algebra then $[abc] : = (a \odot b) \odot c$ induced a
ternary product making ${\mathcal A}$ into a ternary algebra
which will be called trivial. It is known that unital ternary
algebras are trivial and finitely generated ternary algebras are
ternary subalgebras of trivial ternary algebras \cite{B-B-K}.
There are other types of ternary algebras in which one may
consider other versions of associativity (see \cite{MOS1}). Some
examples of ternary algebras are {\it (i)} ``cubic matrices''
introduced by Cayley \cite{CAY} which was in turn generalized by
Kapranov, Gelfand and Zelevinskii \cite{K-G-Z}; {\it (ii)} the
ternary algebra of the polynomials of odd degrees in one variable
equipped with the ternary operation $[p_1p_2p_3]=p_1\odot
p_2\odot p_3$, where $\odot$ denotes the usual multiplication of
polynomials. By a Banach ternary algebra we mean a ternary
algebra equipped with a complete norm $\|.\|$ such that
$\|[abc]\|\leq\|a\| \|b\| \|c\|$.

The stability problem of functional equations originated from a
question of S. Ulam \cite{ULA}, posed in 1940, concerning the
stability of group homomorphisms. In 1941, D. H. Hyers \cite{HYE}
gave a partial affirmative answer to the question of Ulam in the
context of Banach spaces. In 1978, Th. M. Rassias \cite{RAS1}
extended the theorem of Hyers by considering the unbounded Cauchy
difference $\|f(x+y)-f(x)-f(y)\|\leq \varepsilon(\|x\|^ p+\|y\|^
p), ~(\epsilon>0, p\in[0,1))$. The result of Rassias, which is
also true for $p < 0$, has provided a lot of influence in the
development of what we now call {\it Hyers--Ulam--Rassias
stability} of functional equations. In 1992, a generalization of
Th. M. Rassias' theorem, the so-called generalized
Hyers--Ulam--Rassias stability, was obtained by G\u avruta
\cite{GAV} by following the same approach as in \cite{RAS1}.
During the last decades several stability problems of functional
equations have been investigated in the spirit of
Hyers--Ulam--Rassias-G\u avruta. See \cite{CZE, H-I-R, RAS2,
MOS2} and references therein for more detailed information on
stability of functional equations.

Stability of algebraic and topological homomorphisms has been
investigated by many mathematicians, for an extensive account on
the subject see \cite{RAS3}. In \cite{B-M} the authors
investigated the stability of homomorphisms between
$J^*$-algebras associated to the Cauchy equation. Some results on
stability ternary homomorphisms may be found at \cite{A-M, MOS1}.
C. Park \cite{PAR1, PAR2} studied the stability of Poisson
$C^*$-homomorphisms and $JB^*$-homomorphisms associated to the
Jensen equation $2 f(\frac{x+y}{2}) = f(x) + f(y)$ where $f$ is a
mapping between linear spaces. The generalized stability of this
equation was studied by K. Jun and Y. Lee \cite{J-L} and also
\cite{JUN}.

A generalization of the Jensen equation is the equation
\begin{eqnarray*}
r f(\frac{sx+ty}{r}) = s f(x) + t f(y).
\end{eqnarray*}
where $f$ is a mapping between linear spaces and $r, s, t$ are
given constant values (see also \cite{J-M-S}). It is easy to see
that a mapping $f : X \to Y$ between linear spaces with $f(0)=0$
satisfies the generalized Jensen equation if and only if it is
additive; cf. \cite{B-O-P-P}.

In this paper, using some ideas from \cite{B-M, PAR1, MOS1}, we
establish the generalized Hyers--Ulam--Rassias stability of
ternary homomorphisms associated to the generalized Jensen
functional equation. If a ternary algebra $({\mathcal A}, [\;])$
has an identity, i.e. an element $e$ such that $a = [aee] = [eaa]
= [eea]$ for all $a\in {\mathcal A}$, then $a\odot b : = [aeb]$ is
a binary product for which we have
\begin{eqnarray*}
(a\odot b)\odot c = [[aeb]ec] = [ae[bec]] = a\odot (b\odot c)
\end{eqnarray*}
and
\begin{eqnarray*}
a\odot e = [aee] = a = [eea] = e\odot a,
\end{eqnarray*}
for all $a, b, c\in {\mathcal A}$ and so $(A, [\;])$ may be
considered as a (binary) algebra. Conversely, if $(A, \odot)$ is
any (binary) algebra, then $[abc] : = a\odot b\odot c$ makes
${\mathcal A}$ into a ternary algebra with the unit $e$ such that
$a\odot b = [aeb]$. Thus our results may be applied to
investigate of stability of algebra homomorphisms;  see
\cite{PAR3}.

Throughout the paper, ${\mathcal A}$ denotes a ternary algebra,
${\mathcal B}$ is a Banach ternary algebra, $X$ denotes a linear
space and $Y$ represents a Banach space. In addition, we assume
$r, s, t$ to be constant positive integers. If a mapping $f$
satisfies the generalized Jensen equation, then so does
$f(x)-f(0)$. Hence without lose of generality we can assume that
$f(0)=0$.

\section {Main Results }
In this section, we are going to establish the generalized
Hyers--Ulam--Rassias stability of homomorphisms between ternary
algebras associated with the generalized Jensen functional
equation. We start with study of stability of generalized Jensen
equation using different Hyers' sequences from those of
\cite{JUN} and \cite{B-O-P-P}.

\begin{theorem}
Let $f: X\to Y$ be a mapping with $f(0)=0$ for which there exists
a function $\varphi : X\times X \to [0, \infty)$ satisfying
\begin{eqnarray*}
\widetilde{\varphi}(x, y): =\frac{1}{r}\sum_{n=0}^{\infty}
(\frac{r}{s})^{-n} \varphi((\frac{r}{s})^nx,(\frac{r}{s})^ny) <
\infty ,
\end{eqnarray*}
and
\begin{eqnarray}\label{Jen}
\|r f(\frac{sx+ty}{r}) - s f(x) - t f(y)\|\leq\varphi(x,y),
\end{eqnarray}
for all $x, y\in X$. Then there exists a unique additive mapping
$T: X\to Y$ given by $T(x) : =
\lim_{n\to\infty}(\frac{r}{s})^{-n}f((\frac{r}{s})^nx)$ such that
\begin{eqnarray*}
\|f(x) - T(x)\| \leq \widetilde{\varphi}(x,x)
\end{eqnarray*}
for all $x\in X$.
\end{theorem}

\begin{proof} Set $y=0$ in \ref{Jen} to get
\begin{eqnarray*}
\|r f(\frac{sx}{r}) - s f(x)\|\leq\varphi(x,0),
\end{eqnarray*}
whence
\begin{eqnarray*}
\|f(x) - (\frac{r}{s})^{-1} f(\frac{r}{s}x)\|\leq
\frac{1}{r}\varphi(\frac{r}{s}x,0),
\end{eqnarray*}
for all $x\in X$. Assume that for some positive integer $n$
\begin{eqnarray*}
\|f(x) - (\frac{r}{s})^{-n} f((\frac{r}{s})^{n}x)\|\leq
\frac{1}{r}\sum_{k=0}^{n-1}(\frac{r}{s})^{-k}\varphi((\frac{r}{s})^{k+1}x,0),
\end{eqnarray*}
for all $x\in X$. Then
\begin{eqnarray*}
\|f(x) - (\frac{r}{s})^{-n-1} f((\frac{r}{s})^{n+1}x)\| &\leq&
\|f(x) - (\frac{r}{s})^{-n} f((\frac{r}{s})^{n}x)\|\\
&&+\|(\frac{r}{s})^{-n} f((\frac{r}{s})^{n}x) -
(\frac{r}{s})^{-n}(\frac{r}{s})^{-1}f((\frac{r}{s})(\frac{r}{s})^{n}x)
\|\\
&\leq& \|f(x) - (\frac{r}{s})^{-n}
f((\frac{r}{s})^{n}x)\|\\
&&+ (\frac{r}{s})^{-n}\|f((\frac{r}{s})^{n}x) -
 (\frac{r}{s})^{-1}f((\frac{r}{s})(\frac{r}{s})^{n}x)\|\\
&\leq& \frac{1}{r}\sum_{k=0}^{n-1}(\frac{r}{s})^{-k}\varphi((\frac{r}{s})^{k+1}x,0)\\
&&+ (\frac{r}{s})^{-n}\frac{1}{r}\varphi((\frac{r}{s})^{n+1}x,0)\\
&\leq&\frac{1}{r}\sum_{k=0}^{n}((\frac{r}{s})^{-k}\varphi((\frac{r}{s})^{k+1}x,0)
\end{eqnarray*}
for all $x\in X$. Using the induction, we conclude that
\begin{eqnarray}\label{app}
\|f(x) - (\frac{r}{s})^{-n} f((\frac{r}{s})^{n}x)\|\leq
\frac{1}{r}\sum_{k=0}^{n-1}(\frac{r}{s})^{-k}\varphi((\frac{r}{s})^{k+1}x,0),
\end{eqnarray}
for all $x\in X$ and all $n\in N$. Similarly one can show that
\begin{eqnarray*}
\|(\frac{r}{s})^{-n} f((\frac{r}{s})^{n}x) - (\frac{r}{s})^{-m}
f((\frac{r}{s})^{m}x)\|\leq
\frac{1}{r}\sum_{k=m}^{n-1}(\frac{r}{s})^{-k}\varphi((\frac{r}{s})^{k+1}x,0),
\end{eqnarray*}
for all positive integers $n > m$ and all $x\in X$. Hence the
sequence $\{(\frac{r}{s})^{-n} f((\frac{r}{s})^{n}x)\}$ is Cauchy
and so is convergent in the complete space $Y$. Thus we can
define the mapping $T: X\to Y$ by
\begin{eqnarray}\label{lim}
T(x) : = \lim_{n\to\infty}(\frac{r}{s})^{-n}f((\frac{r}{s})^{n}x)
\end{eqnarray}
Replace $x, y$ by $(\frac{r}{s})^{n}x, (\frac{r}{s})^{n}y$,
respectively, in \ref{Jen} to obtain
\begin{eqnarray*}
\|r(\frac{r}{s})^{-n}f(\frac{s(\frac{r}{s})^{n}x-t(\frac{r}{s})^{n}y}{r})
- s (\frac{r}{s})^{-n}f((\frac{r}{s})^{n}x) + t(\frac{r}{s})^{-n}
f((\frac{r}{s})^{n}y)\|\\
\leq(\frac{r}{s})^{-n}\varphi((\frac{r}{s})^{n}x,(\frac{r}{s})^{n}y),
\end{eqnarray*}
for all $x\in X$ and all $n$. Letting $n\to \infty$, we deduce
that $T$ satisfies the generalized Jensen functional equation and
so it is additive. In addition \ref{app} and \ref{lim} yield
\begin{eqnarray*}
\|f(x) - T(x)\| \leq \widetilde{\varphi}(x,x)
\end{eqnarray*}
for all $x\in X$.

We use a standard technique to prove the uniqueness assertion
(see e.g. \cite{MOS2}). First note that for all positive integer
$j$ we have
\begin{eqnarray*}
(\frac{r}{s})^{j}T(x) &=&
(\frac{r}{s})^{j}\lim_{n\to\infty}(\frac{r}{s})^{-n}f((\frac{r}{s})^{n}x)\\
&=&\lim_{n\to\infty}(\frac{r}{s})^{j-n}f((\frac{r}{s})^{n-j} ((\frac{r}{s})^{j}x))\\
&=& T((\frac{r}{s})^{j}x)
\end{eqnarray*}
Now let $T'$ be another additive mapping satisfying $\|f(x) -
T'(x)\| \leq \widetilde{\varphi}(x,x)$ for all $x\in X$. Then
\begin{eqnarray*}
\|T(x) - T'(x)\| &=&
(\frac{r}{s})^{-j}\|T((\frac{r}{s})^{j}x)-T'((\frac{r}{s})^{j}x)\|\\
&\leq&
(\frac{r}{s})^{-j}\|T((\frac{r}{s})^{j}x)-f((\frac{r}{s})^{j}x)\|+
(\frac{r}{s})^{-j}\|f((\frac{r}{s})^{j}x)-T'((\frac{r}{s})^{j}x)\|\\
&\leq& 2
(\frac{r}{s})^{-j}\widetilde{\varphi}((\frac{r}{s})^{j}x,(\frac{r}{s})^{j}x)\\
&=& 2(\frac{r}{s})^{-j}\sum_{k=0}^\infty
(\frac{r}{s})^{-k}\varphi((\frac{r}{s})^{k}
(\frac{r}{s})^{j}x,(\frac{r}{s})^{k}(\frac{r}{s})^{j}x)\\
&=&
2\sum_{k=j}^\infty(\frac{r}{s})^{-k}\varphi((\frac{r}{s})^{k}x,(\frac{r}{s})^{k}x)
\end{eqnarray*}
for all $x\in X$. The right hand side tends to zero as
$j\to\infty$, hence $T(x)=T'(x)$ for all $x\in X$.
\end{proof}

In a similar fashion one may prove the following theorem.
\begin{theorem}
Let $f: X\to Y$ be a mapping with $f(0)=0$ for which there exists
a function $\varphi : X\times X \to [0, \infty)$ satisfying
\begin{eqnarray*}
\widetilde{\varphi}(x, y): =\frac{1}{s}\sum_{n=0}^{\infty}
(\frac{r}{s})^{n}
\varphi((\frac{r}{s})^{-n}x,(\frac{r}{s})^{-n}y) < \infty ,
\end{eqnarray*}
and
\begin{eqnarray*}
\|r f(\frac{sx+ty}{r}) - s f(x) - t f(y)\|\leq\varphi(x,y),
\end{eqnarray*}
for all $x, y\in X$. Then there exists a unique additive mapping
$T: X\to Y$ given by $T(x) : =
\lim_{n\to\infty}(\frac{r}{s})^{n}f((\frac{r}{s})^{-n}x)$ such
that
\begin{eqnarray*}
\|f(x) - T(x)\| \leq \widetilde{\varphi}(x,x)
\end{eqnarray*}
for all $x\in X$.
\end{theorem}

The following proposition gives a sufficient condition in order a
mapping satisfying a Jensen type inequality really to be a
ternary homomorphism.

\begin{proposition} Let $r\neq s$ and $T:{\mathcal A} \to {\mathcal B}$ be a
mapping with $T(\frac{r}{s}x) = \frac{r}{s}T(x), x\in {\mathcal
A}$ for which there exists a function $\varphi : {\mathcal A}^5
\to [0, \infty)$ satisfying
\begin{eqnarray*}
\lim_{n\to\infty}(\frac{r}{s})^{-n}\varphi((\frac{r}{s})^nx,
(\frac{r}{s})^ny, (\frac{r}{s})^n u, (\frac{r}{s})^nv,
(\frac{r}{s})^nw)=0,
\end{eqnarray*}
and
\begin{eqnarray}\label{JenterT}
\|r T(\frac{\mu sx+ \mu ty+[uvw]}{r}) - \mu s T(x) + \mu t
T(y)-[T(u)T(v)T(w)]\|\nonumber\\
\leq\varphi(x,y, u, v, w),
\end{eqnarray}
for all $\mu\in{\mathbb C}, x, y, u, v, w\in{\mathcal A}$.
\end{proposition}
\begin{proof} $T(0)=0$ since $T(0)=\frac{r}{s}T(0)$ and $\frac{r}{s}\neq 1$.
Putting $\mu = 1, u=v=w=0$ and replacing $x, y$ by
$(\frac{r}{s})^nx, (\frac{r}{s})^ny$ in \ref{JenterT}, we get
\begin{eqnarray*}
\|r (\frac{r}{s})^{-n}T((\frac{r}{s})^n\frac{sx+ty}{r}) - s
(\frac{r}{s})^{-n}T((\frac{r}{s})^nx) + t (\frac{r}{s})^{-n}
T((\frac{r}{s})^ny)\|\\
\leq (\frac{r}{s})^{-n}\varphi((\frac{r}{s})^nx,(\frac{r}{s})^ny,
0, 0, 0),
\end{eqnarray*}
Taking the limit as $n\to\infty$ we conclude that $T$ satisfies
the Jensen equation. Hence $T$ is additive. Similarly one can
prove that $T(\mu x)=\mu T(x)$ for all $\mu\in{\mathbb C},
x\in{\mathcal A}$.

Set $x=y=0$ and replace $u, v, w$ by $(\frac{r}{s})^nu,
(\frac{r}{s})^nv, (\frac{r}{s})^nw$ in \ref{JenterT}. Then
\begin{eqnarray*}
\|rT(\frac{[uvw]}{r})-[T(u)T(v)T(w)]\|&=&
(\frac{r}{s})^{-3n}\|r T(\frac{[(\frac{r}{s})^nu((\frac{r}{s})^nv)(\frac{r}{s})^nw]}{r})\\
&&-[T((\frac{r}{s})^nu)T((\frac{r}{s})^nv)T((\frac{r}{s})^nw)]\| \\
&\leq& (\frac{r}{s})^{-3n}\varphi(0, 0, (\frac{r}{s})^nu, (\frac{r}{s})^nv, (\frac{r}{s})^nw)\\
&\leq& (\frac{r}{s})^{-n}\varphi(0, 0, (\frac{r}{s})^nu,
(\frac{r}{s})^nv, (\frac{r}{s})^nw),
\end{eqnarray*}
for all $u, v, w\in {\mathcal A}$. The right hand side tends to
zero as $n\to\infty$, so that $T([uvw]) =
rT(\frac{[uvw]}{r})=[T(u)T(v)T(w)]$ for all $u, v, w\in {\mathcal
A}$. Thus $T$ is a ternary homomorphism.
\end{proof}

\begin{theorem} Let $f:{\mathcal A}\to {\mathcal {\mathcal B}}$ be a
mapping such that $f(0)= 0$, and there exists a function $\varphi
:{\mathcal A}^5 \to [0, \infty)$ such that
\begin{eqnarray}\label{phi}
\widetilde{\varphi}(x, y, u, v, w):
=\frac{1}{r}\sum_{j=0}^{\infty} (\frac{r}{s})^{-j}
\varphi((\frac{r}{s})^j x, (\frac{r}{s})^j y, (\frac{r}{s})^j u,
(\frac{r}{s})^j v, (\frac{r}{s})^j w) < \infty,
\end{eqnarray}
and \begin{eqnarray}\label{Jenter} \|r f(\frac{\mu sx+ \mu
ty+[uvw]}{r}) - \mu s f(x) + \mu t
f(y)-[f(u)f(v)f(w)]\|\nonumber\\
\leq\varphi(x,y, u, v, w),
\end{eqnarray}
for all $\mu\in{\mathbb C}, x, y, u, v, w\in{\mathcal A}$. Then
there exists a unique ternary homomorphism $T:{\mathcal A} \to
{\mathcal B}$ such that
\begin{eqnarray}\label{ft}
\|f(x) - T(x)\|\leq \widetilde{\varphi}(x, x, 0, 0, 0)
\end{eqnarray}
for all $x\in {\mathcal A}$.
\end{theorem}

\begin{proof} Put $u=v=w=0$ and $\mu =1\in{\mathbb T}^1$ in \ref{Jenter}. It follows from theorem
2.1 that there is a unique additive mapping $T:{\mathcal A} \to
{\mathcal B}$ defined by
\begin{eqnarray*}
T(x) = \lim_{n\to\infty}(\frac{r}{s})^{-n} f((\frac{r}{s})^nx),
\end{eqnarray*}
and satisfying \ref{ft} for all $x\in{\mathcal A}$.

Let $\mu\in{\mathbb T}^1$. Replacing $x$ by
$(\frac{r}{s})^{n+1}x$ and $y$ by $0$ in \ref{Jenter}, we get
\begin{eqnarray*}
\|f((\frac{r}{s})^n \mu x)-\mu (\frac{r}{s})^{-1}
f((\frac{r}{s})^{n+1}x)\| \leq
\frac{1}{r}\varphi((\frac{r}{s})^{n+1}x, 0, 0, 0, 0),
\end{eqnarray*}
for all $x\in{\mathcal A}$. It follows from $\|rf(x) -
sf(\frac{r}{s}x)\|\leq \varphi(\frac{r}{s}x,0, 0, 0, 0)$ that
\begin{eqnarray*}
\|f((\frac{r}{s})^n x)- (\frac{r}{s})^{-1}
f((\frac{r}{s})^{n+1}x)\| \leq
\frac{1}{r}\varphi((\frac{r}{s})^{n+1}x, 0, 0, 0, 0),
\end{eqnarray*}
Hence
\begin{eqnarray*}
\|(\frac{r}{s})^{-n}f((\frac{r}{s})^n \mu x)- \mu
(\frac{r}{s})^{-n}f((\frac{r}{s})^{n}x)\|&\leq& (\frac{r}{s})^{-n}
\|f((\frac{r}{s})^n \mu x)- \mu
(\frac{r}{s})^{-1} f((\frac{r}{s})^{n+1}x)\|\\
&&+ (\frac{r}{s})^{-n}\|\mu f((\frac{r}{s})^n x)- \mu
(\frac{r}{s})^{-1}f((\frac{r}{s})^{n+1}x)\|\\
&\leq&
2(\frac{r}{s})^{-n}(\frac{1}{r})\varphi((\frac{r}{s})^{n+1}x, 0,
0, 0, 0),
\end{eqnarray*}
for all $x\in{\mathcal A}$. Taking the limit and using \ref{lim}
and noting that the right hand side tends to zero as $n\to\infty$,
we infer that
\begin{eqnarray*}
T(\mu x) = \lim_{n\to\infty}(\frac{r}{s})^{-n}f((\frac{r}{s})^n
\mu x) = \lim_{n\to\infty} (\mu
(\frac{r}{s})^{-n}f((\frac{r}{s})^nx)) = \mu T(x),
\end{eqnarray*}
for all $x\in{\mathcal A}$. Obviously, $T(0x)=0=0T(x)$.

Next, let $\lambda\in {\mathbb C}$ ($\lambda\neq 0$) and let $M$
be a natural number greater than $4|\lambda|$. Then
$|\frac{\lambda}{M}|<\frac{1}{4}<1-\frac{2}{3}=1/3$. By Theorem 1
of \cite{K-P}, there exist three numbers $\mu_1, \mu_2, \mu_3\in
{\mathbb T}^1$ such that $3\frac{\lambda}{M}=\mu_1+\mu_2+\mu_3$.
By the additivity of $T$ we get $T(\frac{1}{3}x)=\frac{1}{3}T(x)$
for all $ x\in {\mathcal A}$. Therefore,
\begin{eqnarray*}
T(\lambda x)&=& T(\frac{M}{3}\cdot 3 \cdot \frac{\lambda}{M}x)=
MT(\frac{1}{3}\cdot 3\cdot \frac{\lambda}{M}x)
=\frac{M}{3}T(3\cdot \frac{\lambda}{M}x)\\
&=&\frac{M}{3}T(\mu_1x+\mu_2x+\mu_3x)=\frac{M}{3}(T(\mu_1x)+T(\mu_2x)+T(\mu_3x))\\
&=&\frac{M}{3}(\mu_1+\mu_2+\mu_3)T(x)=\frac{M}{3}\cdot 3\cdot \frac{\lambda}{M}\\
&=&\lambda T(x),
\end{eqnarray*}
for all $x\in {\mathcal A}$. So that $T$ is ${\mathbb C}$-linear.

Set $\mu =1$ and $x=y=0$, and replace $u, v, w$ by
$(\frac{r}{s})^nu, (\frac{r}{s})^nv, (\frac{r}{s})^nw$,
respectively, in (\ref{Jenter}). Then
\begin{eqnarray*}
(\frac{r}{s})^{-3n}\|r f((\frac{r}{s})^{3n}\frac{[uvw]}{r})-[f((\frac{r}{s})^{n}u)f((\frac{r}{s})^{n}v)f((\frac{r}{s})^{n}w)]\|\\
\leq (\frac{r}{s})^{-3n}\varphi(0, 0, (\frac{r}{s})^{n}u,
(\frac{r}{s})^{n}v, (\frac{r}{s})^{n}w),
\end{eqnarray*}
for all $u, v, w\in {\mathcal A}$. Then by applying the
continuity of the ternary product $(x,y,z)\mapsto [xyz]$ we deduce
\begin{eqnarray*}
T([uvw])&=& r T(\frac{1}{r}[uvw])\\
&=&\lim_{n\to\infty}(\frac{r}{s})^{-3n} r f((\frac{r}{s})^{3n}\frac{[uvw]}{r})\\
&=&\lim_{n\to\infty}[(\frac{r}{s})^{-n}f((\frac{r}{s})^{n}u)~(\frac{r}{s})^{-n}f((\frac{r}{s})^{n}v)~(\frac{r}{s})^{-n}f((\frac{r}{s})^{n}w)]\\
&=& [T(u)T(v)T(w)],
\end{eqnarray*}
for all $u, v, w\in {\mathcal A}$. Thus $T$ is a ternary
homomorphism satisfying the required inequality.
\end{proof}

\begin{corollary} Let $f:{\mathcal A} \to {\mathcal B}$ be a
mapping such that $f(0) = 0$, and there exist constants $\epsilon
\geq 0$ and $p < 1$ such that
\begin{eqnarray*}
\|r f( \frac{\mu s x + \mu t y+[uvw]}{r})- \mu s f(x) -
\mu t f(y)-[f(u)f(v)f(w)]\|\\
\leq \epsilon (\|x\|^p + \|y\|^p + \|u\|^p+ \|v\|^p + \|w\|^p),
\end{eqnarray*}
for all $\mu\in{\mathbb T}^1$, and all $x, y\in{\mathcal A}$ and
all $u, v, w \in{\mathcal A}$. Then there exists a unique ternary
homomorphism $T:{\mathcal A} \to {\mathcal B}$ such that
\begin{eqnarray*}
\|f(x) - T(x)\|\leq
\frac{2r^{-p}\epsilon\|x\|^p}{r^{1-p}-s^{1-p}},
\end{eqnarray*}
for all $x \in{\mathcal A}$.
\end{corollary}

\begin{proof} Define $\varphi(x, y, u, v, w) = \epsilon (\|x\|^p +
\|y\|^p + \|u\|^p + \|v\|^p +\|w\|^p)$, and apply Theorem 2.4.
\end{proof}

\begin{remark}
When $p>1$, one may use the same techniques used in the proof of
Theorem 2.4 (Theorem 2.5) to get a result similar to Corollary
2.5.
\end{remark}

The following corollary may applied in the case that our ternary
algebra is linearly generated by ``idempotents'', i.e. elements
$u$ with $u^3 = u$.
\begin{corollary}
Let ${\mathcal A}$ be linearly spanned by a set $S\subseteq
{\mathcal A}$ and $f:{\mathcal A} \to {\mathcal B}$ be a mapping
with $f(0)=0$ satisfying $f((\frac{r}{s})^{2n}[s_1s_2z]) =
[f((\frac{r}{s})^{n}s_1)f((\frac{r}{s})^{n}s_2)f(z)]$ for all
sufficiently large positive integers $n$, and all $s_1,s_2\in S,
z\in{\mathcal A}$. Suppose that there exists a function $\varphi$
fulfilling (\ref{phi}), and
\begin{eqnarray*}
\|r T(\frac{\mu sx+ \mu ty+[uvw]}{r}) - \mu s T(x) - \mu t
T(y)-[T(u)T(v)T(w)]\|\leq\varphi(x, y, u, v, w),
\end{eqnarray*}
for all $\mu\in{\mathbb T}^1$ and all $x, y\in{\mathcal A}$. Then
there exists a unique ternary homomorphism $T:{\mathcal A} \to
{\mathcal B}$ satisfying \ref{ft} for all $x\in{\mathcal A}$.
\end{corollary}
\begin{proof}
By the same arguing as in the proof of Theorem 2.4, there exists a
unique linear mapping $T:{\mathcal A} \to {\mathcal B}$ given by
\begin{eqnarray*}
T(x) : = \lim_{n\to\infty}(\frac{r}{s})^{-n}f((\frac{r}{s})^{n}x),
~~~~~~ (x\in{\mathcal A}).
\end{eqnarray*}
such that
\begin{eqnarray*}
\|f(x) - T(x)\|\leq \widetilde{\varphi}(x, x, 0, 0, 0),
\end{eqnarray*}
for all $x\in{\mathcal A}$. We have
\begin{eqnarray*}
T([s_1s_2z]) &=& \lim_{n\to\infty}(\frac{r}{s})^{-2n} f([((\frac{r}{s})^{n}s_1)((\frac{r}{s})^{n}s_2)z])\\
&=& \lim_{n\to\infty} [(\frac{r}{s})^{-n}f((\frac{r}{s})^{n}s_1) (\frac{r}{s})^{-n}f((\frac{r}{s})^{n}s_2)f(z)])\\
&=& [T(s_1)T(s_2)f(z)]
\end{eqnarray*}
By the linearity of $T$ we have $T([xyz]) = [T(x)T(y)f(z)]$ for
all $x, y, z\in {\mathcal A}$. Therefore
$(\frac{r}{s})^{n}T([xyz])= T([xy((\frac{r}{s})^{n}z)]) =
[T(x)T(y)f((\frac{r}{s})^{n}z)]$, and so
\begin{eqnarray*}
T[xyz])&=&
\lim_{n\to\infty}(\frac{r}{s})^{-n}[T(x)T(y)f((\frac{r}{s})^{n}z)]\\
&=&[T(x)T(y)\lim_{n\to\infty}(\frac{r}{s})^{-n}f((\frac{r}{s})^{n}z)]\\
&=& [T(x)T(y)T(z)],
\end{eqnarray*}
for all $x,y,z\in{\mathcal A}$.
\end{proof}

\begin{theorem} Suppose that $f:{\mathcal A} \to {\mathcal B}$ is
a mapping with $f(0)=0$ for which there exists a function
$\varphi: A^{5}\to [0, \infty)$ fulfilling (\ref{phi}), and
(\ref{Jenter}) holds for $\mu=1, {\bf i}$ and all $x\in {\mathcal
A}$. Then there exists a unique ternary homomorphism $T:{\mathcal
A} \to {\mathcal B}$ such that
\begin{eqnarray*}
\|f(x)-T(x)\|\leq \widetilde{\varphi}(x, x, 0, 0, 0),
\end{eqnarray*}
for all $x\in{\mathcal A}$.
\end{theorem}
\begin{proof} Put $u=v=w=0$ and $\mu=1$ in (\ref{Jenter}). By the same arguing as Theorem 2.4
we infer that there exists a unique additive mapping $T:{\mathcal
A} \to {\mathcal B}$ given by
\begin{eqnarray*}
T(x) : = \lim_{n\to\infty}(\frac{r}{s})^{-n}f((\frac{r}{s})^{n}x),
\end{eqnarray*}
and satisfying \ref{ft} for all $ x\in {\mathcal A}$. By the same
reasoning as in the proof of the main theorem of \cite{RAS1}, the
mapping $T$ is ${\mathbb R}$-linear.

Replace $x$ by $(\frac{r}{s})^{n}x$, $y$ by $0$, and put $\mu={\bf
i}$ and $u=v=w=0$ in (\ref{Jenter}). Then
\begin{eqnarray*}
(\frac{r}{s})^{-n}\|\frac{r}{s}f({\bf i}(\frac{r}{s})^{n-1}x)
-{\bf i} f((\frac{r}{s})^{n}x)\| \leq
\frac{1}{s}(\frac{r}{s})^{-n}\varphi((\frac{r}{s})^{n}x, 0, 0, 0,
0),
\end{eqnarray*}
for all $x\in{\mathcal A}$. The right hand side tends to zero as
$n\to\infty$, hence
\begin{eqnarray*}
T({\bf
i}x)=\lim_{n\to\infty}(\frac{r}{s})^{-n+1}f((\frac{r}{s})^{n-1}{\bf
i}x)=\lim_{n\to\infty}{\bf
i}(\frac{r}{s})^{-n}f((\frac{r}{s})^{n}x)={\bf i}T(x),
\end{eqnarray*}
for all $ x\in {\mathcal A}$.

For every $\lambda\in {\mathbb C}$ we can write
$\lambda=\alpha_1+{\bf i}\alpha_2$ in which
$\alpha_1,\alpha_2\in{\mathbb R}$. Therefore
\begin{eqnarray*}
T(\lambda x)&=&T(\alpha_1x+{\bf i}\alpha_2x)=\alpha_1T(x)+\alpha_2T({\bf i}x)\\
&=&\alpha_1T(X)+{\bf i}\alpha_2T(x)=(\alpha_1+{\bf i}\alpha_2)T(x)\\
&=&\lambda T(x),
\end{eqnarray*}
for all $x\in {\mathcal A}$. Thus $T$ is ${\mathbb C}$-linear.
\end{proof}
\begin{remark}
There are similar results by considering $r, t$ instead of $r, s$
throughout the paper (especially in Theorems 2.1 and 2.2). In the
case $r=s=t=1$, the generalized Jensen equation turns out the
Cauchy equation.
\end{remark}


\begin{thebibliography}{10}
\bibitem{A-M} M. Amyari and M. S. Moslehian, \textit{Approxiamte homomorphisms of ternary
semigroups}, to appear in Lett. Math. Phys.,
http://arxiv.org/math-ph/0511039.
\bibitem{B-M} C. Baak and M. S. Moslehian, \textit{Stability of
$J^*$-homomorphisms}, Nonlinear Analysis (TAM), 63 (2005), 42-48.
\bibitem {B-B-K} N. Bazunova, A. Borowiec and R. Kerner, \textit{Universal differential calculus on ternary algebras}, Lett. Math. Phys. 67 (2004), no. 3, 195--206.
\bibitem {B-O-P-P} D.-H. Boo, S.-Q. Oh, C.-G. Park and J.-M. Park, \textit{Generalized Jensen's equations in Banach modules over a $C^*$-algebra and its unitary group}, Taiwanese J. Math. 7 (2003),
no. 4, 641--655.
\bibitem{CAY} A. Cayley , Cambridge Math. Journ. {\bf 4}, p. 1 (1845).
\bibitem{CZE} S. Czerwik, \textit{Stability of Functional Equations of Ulam--Hyers--Rassias Type}, Hadronic Press, 2003.
\bibitem{GAV} P. G\u avruta, \textit{A generalization of the Hyers--Ulam--Rassias stability of approximately additive mappings},
 J. Math. Anal. Appl. 184 (1994), 431--436.
\bibitem{HYE} D. H. Hyers, \textit{On the stability of the linear functional equation}, Proc. Nat. Acad. Sci. U.S.A. 27 (1941), 222--224.
\bibitem{H-I-R} D. H. Hyers, G. Isac and Th. M. Rassias, \textit{Stability of Functional Equations in Several Variables}, Birkh\"auser, Basel, 1998.
\bibitem{K-P} R. V. Kadison and G. K. Pedersen, \textit{Means and convex combinations of unitary operators}, Math. Scan. 57 (1985), 249--266.
\bibitem{K-G-Z} M. Kapranov, I. M. Gelfand, and A. Zelevinskii,\textit{Discriminants, Resultants and Multidimensional Determinants},
Birkh\"auser, Berlin, 1994.
\bibitem{J-L} K. Jun and Y. Lee, \textit{A generalization of the
Hyers--Ulam--Rassias stability of Jensen's equation}, J. Math.
Anal.  Appl. 238(1999), 305--315.
\bibitem {JUN} S.-M. Jung, \textit{Hyers-Ulam-Rassias stability of Jensen's equation and its application}, Proc. Amer. Math. Soc. 126 (1998), no. 11, 3137--3143.
\bibitem{J-M-S} S.-M. Jung, M. S. Moslehian, P. K. Sahoo, \textit{Satbility of generalized Jensen equation on restricted
domains}, preprint, http://arxiv.org/math.FA/0511320.
\bibitem{MOS1} M. S. Moslehian, \textit{Approximate $C^*$-ternary ring homomorphisms associated to the Trif
equation}, preprint, http://arxiv.org/math.FA/0511539.
\bibitem{MOS2} M. S. Moslehian, \textit{Approximately vanishing of topological cohomology groups}, J. Math. Anal. Appl., in press, http://arxiv.org/math.FA/0501015
\bibitem {PAR1}  C.-G. Park, \textit{Homomorphisms between Poisson $JC^*$-algebras} Bull. Braz. Math. Soc. (N.S.) 36 (2005), no. 1, 79--97.
\bibitem {PAR2}  C.-G. Park, \textit{Approximate homomorphisms on $JB^*$-triples}, J. Math. Anal. Appl. 306 (2005), no. 1, 375–-381.
\bibitem {PAR3}  C.-G. Park, \textit{On an approximate automorphim on a $C^*$-algebra}, to appera in Proc. Amer. Math.
Soc.
\bibitem{RAS1} Th. M. Rassias, \textit{On the stability of the linear mapping in Banach
spaces}, Proc. Amer. Math. Soc. 72 (1978), 297--300.
\bibitem{RAS2} Th. M. Rassias (ed.), \textit{Functional Equations, Inequalities and Applications}, Kluwer Academic Publishers, Dordrecht, Boston, Londonm, 2003.
\bibitem{RAS3} Th. M. Rassias, \textit{Approximate homomorphisms}, Aequationes Math. 44 (1992), no. 2-3, 125--153.
\bibitem{ULA} S. M. Ulam, \textit{Problems in Modern Mathematics}, Chapter VI, Science Editions, Wiley, New York, 1964.
\end{thebibliography}
\end{document}